\magnification=\magstep1

\font\large=cmbx10 scaled \magstep1
\pageno=1

\normalbaselines 
\overfullrule=0pt           
\footline{Appeared in {\it R. C. P. 25 }, Pr\'epubl.\ IRMA,
Strasbourg, vol.\ 47 (1995) 269-274. \hfil\folio\hfil }

\def\ie{{\it i.e.\ \/}}
\def\cf{{\it cf.\ \/}}

\vskip2cm

\noindent
\centerline{{\large A proof via the Seiberg-Witten moduli space of
Donaldson's }}
\medskip\noindent
\centerline{{\large theorem on smooth 4-manifolds
with definite intersection forms}}
\bigskip\noindent
\vskip1cm
\centerline{Mikhail Katz}
\bigskip\noindent
\vskip2cm

Most of what follows was explained to me by D. Kotschick,
with additional clarifications by T. Delzant, J.-C. Sikorav,
and K. Wojciechowski.  
In [1], the existence of such a proof is attributed to P. Kronheimer 
and others.
Compared to the original proof of Donaldson's theorem, the proof
using the new moduli space is essentially trivial,
which is what motivated a nonspecialist to present this
exposition for nonspecialists in the field.
\medskip\noindent
{\bf Theorem } (S. Donaldson [9]).  Let $X$ be a smooth oriented
4-manifold.  Suppose that the intersection form of $X$ is negative
definite.  Then it is minus the identity.
\medskip\noindent
We perform a surgery, without changing the intersection form in
2-dimensional homology modulo torsion, to reduce to the case of $X^4$
with $b_1=0$ (for details, see section 6).

Consider the Seiberg-Witten moduli space associated with a complex
line bundle $L$ satisfying $c_1(L)=w_2(TX)$ (mod 2) (\cf [3], [6]).
The Seiberg-Witten equations are
$$D_A\phi=0,\ 	F_A^+=i\sigma(\phi,\phi)
\eqno{(1)}
$$
in unknowns $(A,\phi)$.  Here $\phi\in\Gamma(V_{+})$
is a positive spinor, $V_{+}$ is a spin$^{c}$
structure of determinant $L$,
while $D_A$ is the Dirac operator built from a connection $A$
on $L$ and the Levi-Civita connection on $TX$, and $\sigma$
is a quadratic form.
The self-dual part $F_A^{+}$ of the curvature of $A$ is 
in $\Lambda^{2,+}(X)$.
By Clifford multiplication, $\Lambda^{2,+}$ acts by endomorphisms
of the bundle $V_{+}$.  At a point, such an 
endomorphism is given by a 2 by 2 matrix.  To build such a 
matrix out of $\phi$, think of $\phi$ as a column vector with
two components.
Then $\sigma(\phi,\phi)$
is like $(\phi \phi^*)_{0}$, the trace-free part of
the 2 by 2 matrix obtained as
the product of $\phi$ by its conjugate transpose.

The gauge group $G= Map(X, U(1))$ acts on solutions by 
$g.A=A-2d\log g,\ g.\phi=g\phi$ (multiplication by
complex scalars).
The Seiberg-Witten moduli space is 
$$M_L=\hbox{solutions}/G.$$
The action of $G$ is free except at ``reducible points''
$\phi=0$, \ie $(A,\phi)=(A,0)$.
Choose a basepoint $x_0\in X$.
The based gauge group $G_0$ is the subgroup of $G$
 defined as the set of $g\in G$
such that $g(x_0)=1\in U(1)$
(gauge transformations fixing the fiber over the basepoint).
We have an exact sequence
$G_0 \rightarrow G \rightarrow U(1).$
The action of $G_0$ is free.
The based moduli space $M_0$ is the quotient
of the space of solutions of SW equations by $G_0$.

To prove Donaldson's theorem, we argue by contradiction.
Suppose the negative definite intersection form is not minus
the identity.  The argument is in 5 steps:
\medskip
1. We specify an $L$ defining a moduli space $M_L$
of positive virtual dimension (using

\hskip.1in
the fact that the intersection form is not minus the identity).

2. The reducible point is unique (we use $b_1=0$ here).

3. We truncate around the reducible point to arrive at a 
contradiction with a

\hskip.1in
standard result on characteristic numbers
called Pontrjagin's theorem.

4. We perturb the second equation to ensure genericity,
and verify the existence and 

\hskip.1in
uniqueness of the reducible point
for the  perturbed equations.

5. We perturb the first equation to ensure the smoothness
of the based moduli 

\hskip.1in
space $M_0$ at the reducible point.

\medskip
The arithmetic source of Donaldson's theorem
is a remark of N. Elkies [1]:
\medskip\noindent
{\bf Theorem.}  The identity is the only bilinear unimodular 
positive definite form ( , ) over $\bf Z$ which does not admit a 
vector $w\in {\bf Z}^n$ satisfying the following 2 properties:

(a) $(w,w) < n$;

(b) for all $v\in Z^n$ one has $(v,v+w)=0$ (mod 2).
\medskip\noindent
Such a $w$ will be called a short characteristic vector.
Now
$(w,w)$ is congruent to the signature modulo 8 (\cf [2]).
Thus in the positive definite case, any non-diagonal form
admits a $w$ such that
$$(w,w)= n-8k {\rm \ with\ } k \geq 1.
\eqno{(2)}
$$
\bigskip\noindent
{\large 1. Choice of $L$ defining
a moduli space of positive dimension}  
\medskip\noindent
The condition $c_1(L)=w_2(TX)$ of the existence of a 
spin$^c$ structure means (by Wu's theorem)
that $c_1(L)$ is a characteristic
vector of the intersection form of $X$.  We choose the short
one, or more precisely any class whose reduction modulo
torsion is the short vector.  Then (2) gives
$$c_1(L)^2 = -|c_1(L)^2| = -b_2 + 8k,\ k\geq 1
\eqno{(3)}
$$
since the intersection form is negative definite.
The virtual dimension of the SW moduli space is 
(\cf [6])
one quarter of
$$c_1(L)^2-(2\chi+3\sigma)=-b_2+8k-(4-4b_1-b_2)
= 8k-4+4b_1 > 0 {\rm \ if\ } k\geq 1.
\eqno{(4)}
$$
\bigskip\noindent
{\large 2. Uniqueness of the reducible point}  
\medskip\noindent
Reducible solutions $(A,0)$ of equations (1)
are characterized by $F_A^+=0$ \ie
$*F_A=-F_A$.  Since the curvature form $F_A$ is closed, applying $d$ 
we see that $F_A$ is harmonic.  
Existence is immediate since every harmonic form is anti-self-dual
by hypothesis (\cf section 4 below for more details).  By Gauss-Bonnet
the cohomology class $[F_A]=2i\pi c_1(L)$ is prescribed by the choice of $L$,
hence $F_A$ is unique.  Now suppose there are 2 connections $A$ and $A'$
with the same curvature form.  Their difference is therefore
a closed 1-form, hence exact ($b_1=0$).  Thus 
$$A - A' = i df
\eqno{(5)}
$$
and the gauge transformation $g(x)=e^{if(x)/2}$ establishes
the gauge equivalence of $A$ and $A'$.

For example, the flat connection on
a line bundle $L$ whose  $c_1(L)$ is torsion, is unique.
The flat connections up to gauge equivalence correspond to
representations of the fundamental group in $U(1)= R/Z.$
If $b_1=0$, different representations of $\pi_1(X)$ define 
different line bundles, and hence the flat connection on $L$ 
is unique.  Thus for the Enriques surface,  $\pi_1=Z/2Z$, 
there are 2 representations in $U(1)$ hence
two flat connections, but the nontrivial one lives 
on the canonical bundle.  The latter is nontrivial since the
surface is not spin.
\bigskip\noindent
{\large 3.  Truncating around the reducible point
and contradiction with Pontrjagin's theorem}
\medskip\noindent
Consider again
the moduli space $M_L$ of dimension $2k-1$ from formula (4).
The compactness of $M_L$ is established using a Weitzenbock
formula and a $C^0$ estimate on the size of $\phi$ (\cf [6]).
Assume that away from the reducible point,
$M_L$ is nonempty and smooth (see section 4).

Consider the based moduli space $M_0$ which is the quotient of the
space of solutions of SW equations by $G_0$, gauge transformations
fixing the fiber over a basepoint.  Note that dim $M_0 = 2k$.  Let
$p\in M_0$ be the preimage of the reducible point $(A,0) \in M$.  The
complement of a small neighborhood of $p$ is then a manifold (with
boundary), which is compact since the moduli space is compact.
Its existence will lead to a contradiction.

Assume $M_0$ is smooth at $p$ (see section 5).
Choose a metric
on $M_0$ invariant under the action of $U(1)$.  
The induced linear action in the tangent space
$T_p M_0$ at $p$
is free, for otherwise some vector would have a nontrivial
finite stabilizer.  
Via the exponential map this would produce a point
in $M_0\setminus\{p\}$ with a nontrivial 
stabilizer, contradicting
the freeness of the action of all of G on irreducible solutions.
Multiplication by $i\in U(1)$ thus defines a complex structure
on $T_p M_0$.  
The action of $U(1)$ in the tangent space is scalar.
Factoring by the $U(1)$ action we obtain the
standard quotient 
$$S^{2k-1}/S^1=CP^{k-1}.
\eqno{(6)}
$$
The exponential map at $p$ is equivariant with respect to the $U(1)$
action.  Therefore the quotient of a small distance sphere centered at
$p$ by $U(1)$ is still $CP^{k-1}.$ We delete from $M$ the
neighborhood of the reducible point bounded by the $CP^{k-1}$,
to obtain a $(2k-1)$-dimensional manifold $V$ whose boundary is
$CP^{k-1}$.  Note that $V$ is smooth since the deleted neighborhood
contains the only singular point of $M$.  For example, if $k=1$ we
obtain a compact 1-dimensional manifold whose boundary is a single
point, which is already a contradiction.

Consider the
circle bundle over $V$ defined by the projection $M_0 \rightarrow M$.
Its restriction to $CP^{k-1}$ is the Hopf fibration, of non-zero
second Stiefel-Whitney class $w_2$.
Hence its {\it number} is nonzero: $w_2^{k-1}[CP^{k-1}]\not=0$.
But by Pontrjagin's theorem, all such numbers have to vanish,
as the fibration extends over all of $V$ (\cf [4], p.\ 52;
the argument given here for the tangent bundle works also
for the Hopf fibration).  Note
that we have made no use of the orientability of $V$.
The contradiction proves that a non-diagonal intersection form
on a smooth 4-manifold could not have existed in the first place.
\bigskip\noindent
{\large 4.
Existence and uniqueness of the reducible point
for the  perturbed equations}
\medskip\noindent
In [6] it is shown that the perturbed SW equations
$$D_A\phi=0,\ 	F_A^+ - i \sigma(\phi,\phi) = e
\eqno{(7)}
$$
for generic $e\in \Lambda^+$, have a smooth moduli space 
of the dimension predicted by the index theorem,
using the surjectivity of the linearized operator
and the existence of a suitable slice for the action of $G$.
Here one needs
the unique continuation property for spinors in the
kernel of $D_A$ (\cf [7]).

We now check that it always contains a reducible point,
i.e. that the perturbed equation 
$$F_A^+ = e, {\rm where\ } e\in \Lambda^+,
\eqno{(8)}
$$
has a solution.  Then if $M_0$ is smooth at this point
(see section 5), we can conclude that the (irreducible)
moduli space is non-empty.
Simultaneously we check that the reducible point is unique, 
so that step 2 above goes through when the
equations are perturbed.

Consider the Hodge decomposition
$$e= df + *dg + h
\eqno{(9)}
$$
where $h$ is a harmonic 2-form, and $f$ and $g$ are 1-forms
unique up to adding exact 1-forms (since $b_1=0$).
In our set-up a harmonic self-dual form is necessarily 0
hence $h=0$ and the equation $*e=e$ implies $df=dg.$
Thus $e=(1+*)df=d^+(f)$.  Now pick any connection $A_0$ and find
$f$ such that 
$$d^+(f)=e - F_{A_0}^+
\eqno{(10)}
$$
from the Hodge decomposition of the right hand side.
The connection $A_0 + f$  then solves the perturbed equation.
The solution is unique up to adding an exact 1-form, i.e.
up to gauge transformation.  
\bigskip\noindent
{\large 5. Smoothness of the based moduli space at the reducible
point}
\medskip\noindent
The linearisation of the equations at $(A,0)$ is 
$D_A:\Gamma(V_{+})\rightarrow \Gamma(V_{-})$.
This operator is not always surjective as it is possible
to have `harmonic' spinors of both chiralities.
We perturb the first equation by adding a 1-form $c$ to the
operator:
$$(D_A+c)\phi=0,\ 	F_A^+ - i\sigma(\phi,\phi)=e.
\eqno{(11)}
$$
The equations are still gauge-invariant.
We need to verify that the moduli space is still compact.
But the operator $D_A+c=D_{A'}$ with $A'=A+2c$ is still
of the same type.  Applying the Weitzenbock formula to
$D_{A'}$ as in [6], we obtain the necessary $C^0$ estimate
for $\phi$, containing an additional term $|2d^{+}c+e|$ besides
the scalar curvature as follows.  At a point where $|\phi |$ is
a maximum, we have as in [6],
$$\matrix{0&\leq&\Delta|\phi|^2\leq
-{s\over 2}|\phi |^2 + <F_{A'}^{+}\phi,\phi>
\hfill
\cr 
 &=&
-{s\over 2}|\phi |^2 + <(F_{A}^{+}+2d^{+}c)\phi,\phi>
\hfill
\cr 
 &=&
-{s\over 2}|\phi |^2 +<i\sigma(\phi,\phi)\phi,\phi>+
<(2d^{+}c+e)\phi,\phi>
\hfill
\cr 
 &\leq&
-{s\over 2}|\phi |^2 -{1\over 4}|\phi|^4+|2d^{+}c+e|\ |\phi|^2\ .
\hfill
\cr}$$
To choose a suitable 1-form $c$, we suppose for simplicity that
ind$(D_A)=0$
and the kernel is 1-dimensional.
Let $\alpha\in Ker(D_A)$ and 
$\beta\in Im(D_A)^{\perp}\subset\Gamma(V_{-})$.
By the unique continuation property (\cf [7]), there exists
a point $x$ such that $\alpha(x)\not=0$ and $\beta(x)\not=0$.
We choose $c$ so that $c(x).\alpha(x)=\beta(x)$ for Clifford
multiplication.
Choose a function $\psi$ with support near $x$, and let
$\psi_\epsilon=\epsilon\psi$.  Then
$D_A+\psi_\epsilon c$ is invertible
for small $\epsilon$ (\cf [8]).
\bigskip\noindent
{\large 6. Remarks}
\medskip\noindent

1.  If one changes the orientation of $X$ so that the intersection
form is positive definite (while keeping $F^+$ in the SW equations),
the calculation changes.  There is no difficulty in producing
a positive dimensional moduli space, but there is no reducible
point since there are no anti-self-dual harmonic forms.
\medskip
2.  Surgery along a loop representing the free
part of $H_1(X)$ does not change the intersection form
in $H_2(X,Z)$
modulo torsion.
Indeed, let $S^1$ be a loop in $X$ representing a nonzero
class in $H_1(X,R)$.  Let $X_{\_}$ be a neighborhood of $S^1$
(homeomorphic to $S^1\times B^3$)
and $X_{+}$ its complement.  The Mayer-Vietoris sequence gives
$$0 \rightarrow H_2(X_{\_})\oplus H_2(X_{+})\rightarrow 
H_2(X) \rightarrow 0
\eqno{(12)}
$$
where the last arrow is zero by assumption on $S^1$,
and the first arrow replaces a homomorphism which
is zero because $H_3(X) \rightarrow H_2(S^1\times S^2)$ is surjective
(it suffices to consider the ``dual" surface transverse
to $S^1$).  After the surgery replacing $X_{\_}$ by 
$B^2\times S^2$, 
we obtain the exact sequence
$$
0 \rightarrow H_2(S^1\times S^2) \rightarrow 
H_2(B^2\times S^2)\oplus H_2(X_+) 
\rightarrow H_2(X') \rightarrow 0
\eqno{(13)}
$$
where $X'$ is the result of applying the standard surgery on $X$
along $S^1$.  Here the second arrow is injective because
$H_2(S^1\times S^2)={\bf R}$ is isomorphic to $H_2(B^2\times S^2).$

In passing from 0 to $H_2(S^1\times S^2)$ we have increased the
rank by 1, and in passing from $H_2(X_{\_})$ to $H_2(B^2\times S^2)$
we have increased the rank by 1.  Since the Euler
characteristic of the sequence is still 0, it follows that
the rank of $H_2(X)$ is unchanged.  Moreover, neither is the
intersection form 
(in $H_2(X,Z)$ modulo torsion), since we can perform the
surgery so as to avoid a family of 2-cycles representing the generators
of $H_2(X,{\bf Z})$, and in particular leave their intersections unchanged.
\bigskip
3. Can we avoid the equations (11) and the proof of the
smoothness of the based moduli space at the reducible point?
At the reducible point,
the linearisation of the perturbed SW equations (7)
is $D_A(\phi)$.
Perturbing the metric does not solve the problem here,
for it is not known whether $D_A$ is surjective for a generic
metric (\cf [5]).
If $D_A$ is not surjective, then the Kuranishi model (I am now using
the notation from 
[6]) for the moduli
space near the reducible point is $\psi^{-1}(0)/U(1)$ where
$$\psi : Ker(D_A)\rightarrow Coker(D_A) 
\eqno{(14)}
$$
and the $U(1)$ action on the spinors is scalar.

Here the advantage of the Kuranishi model
seems to be that we can still take a sphere of radius $\epsilon$
in $Ker(D_A)$, which will be transverse to $\psi^{-1}(0)$ for almost
all $\epsilon$.  The quotient by $U(1)$ will then give a smooth manifold 
(the analogue of the projective space in formula (6) of section 3).
Consider the restriction of the $U(1)$ bundle over the irreducible
moduli space to this manifold.   Does
this bundle have nontrivial Stiefel-Whitney numbers?
\bigskip
4. The doubly perturbed equations (11) can be avoided by incorporating
the term $2d^+c$ into the right hand side of the second equation of (7).
Since the purpose of the perturbation is to make the {\it first}
equation surjective, that's the equation we chose to perturb.

\bigskip\noindent
{\large References}
\medskip\noindent
[1] N. Elkies, A characterization of the ${\bf Z}^n$ lattice,
Math.\ Research Letters 2 (1995).
\medskip\noindent
[2] J.-P. Serre, A Course in Arithmetic.
\medskip\noindent
[3] E. Witten, Monopoles and four-manifolds, Math.\ Research Letters
1 (1994) 769-796.
\medskip\noindent
[4] J. Milnor and J. Stasheff, Characteristic classes.
Princeton University Press, 1974.
\medskip\noindent
[5] D. Kotschick, 
{Non--trivial harmonic spinors on generic algebraic surfaces},
Proc.\ of the AMS.
\medskip\noindent
[6] P. Kronheimer and T. Mrowka, The genus of embedded surfaces
in the projective plane, Math. Research Letters 1 (1994) 797-808.
\medskip\noindent
[7] B. Booss-Bavnbek and K. Wojciechowski, Elliptic boundary problems for
Dirac operators. Birkhauser, 1993.
\medskip\noindent
[8] K. Wojciechowski and S. Klimek, Unique continuation property
and surjectivity of elliptic operators of order 1, preprint.
\medskip\noindent
[9] S. Donaldson, The orientation of Yang-Mills moduli spaces
and 4-manifold topology, J. Differential Geometry 26 (1987) 397-428.
\bigskip
\vfill
\noindent


\noindent
\vfill\eject\end